\documentclass[11pt]{article}
\usepackage{amsmath,amsfonts,amssymb}
\newtheorem{theorem}{Theorem}[section]

\newtheorem{lemma}[theorem]{Lemma}

\newtheorem{definition}[theorem]{Definition}
\numberwithin{equation}{section}
\let\G=\Gamma
\let\L=\Lambda

\let\g=\gamma

\let\S=\Sigma
\let\O=\Omega

\let\l=\lambda
\let\L=\Lambda
\begin{document}
\title{Analytic Continuation of Holomorphic Correspondences and
  Equivalence of Domains in ${\mathbb C}^n$.} 
\author{Rasul Shafikov \\{\small \tt shafikov@math.sunysb.edu}}
\date{\today}
\maketitle

\begin{abstract}
The following result is proved: Let $D$ and $D'$ be bounded domains in
$\mathbb C^n$, $\partial D$ is smooth, real-analytic, simply
connected, and $\partial D'$ is connected, smooth, real-algebraic. Then there
exists a proper holomorphic correspondence $f:D\to D'$ if and only if
$\partial D$ and $\partial D'$ a locally CR-equivalent. This leads
to a characterization of the equivalence relationship between bounded domains in
$\mathbb C^n$ modulo proper holomorphic correspondences in terms of
CR-equivalence of their boundaries. 
\end{abstract}

\section{Definitions and Main Results.}

Following Stein \cite{st1} we say that a holomorphic correspondence
between two domains
$D$ and $D'$ in ${\mathbb C}^n$ is a complex-analytic set $A \subset
D\times D'$ which satisfies: (i) $A$ is of pure complex dimension $n$,
and (ii) the natural projection $\pi:A\to D$ is proper. A
correspondence $A$ is called proper, if in addition the natural projection
$\pi': A\to D'$ is proper. We may also think of $A$ as the graph of
a multiple valued mapping defined by $F:=\pi'\circ\pi^{-1}$. 
Holomorphic correspondences were studied for instance, in \cite{bb}, 
\cite{bb2}, \cite{be2}, \cite{bsu}, \cite{dp1}, \cite{dfy}, \cite{p3},
and \cite{v}. 

In this paper we address the following question: Given two domains $D$
and $D'$, when does there exist a proper holomorphic correspondence 
$F:D\to D'$? 
Note (see \cite{st2}) that the existence of a correspondence $F$ 
defines an equivalence relation $D \sim D'$. This equivalence relation
is a natural generalization of biholomorphic equivalence of domains in
${\mathbb C}^n$.

To illustrate the concept of equivalence modulo holomorphic
correspondences, consider domains of the form 
$\Omega_{p,q}= \{|z_1|^{2p}+|z_2|^{2q}<1 \}, \ \ p,q\in{\mathbb Z}^+$.
Then $f(z)=({z_1}^{p/s},{z_2}^{q/t})$ is a proper holomorphic 
correspondence between $\Omega_{p,q}$ and $\Omega_{s,t}$, while a 
proper holomorphic map from  $\Omega_{p,q}$ to $\Omega_{s,t}$ exists
only if $s|p$ and $t|q$, or $s|q$ and $t|p$. For details see \cite{bd}
or \cite{l}. 

The main result of this paper is the following theorem.

\begin{theorem}\label{t1}
Let $D$ and $D'$ be bounded domains in ${\mathbb C}^n$, $n>1$. Let 
$\partial D$ be smooth, real-analytic, connected and simply connected,
and let $\partial D'$ be connected, smooth, real-algebraic. Then there
exists a proper holomorphic correspondence $F:D\to D'$ if and only if
there exist points $p\in\partial D$, $p'\in\partial D'$, neighborhoods
$U\ni p$ and $U'\ni p'$, and a biholomorphic map $f:U\to U'$, such
that $f(U\cap\partial D)= U'\cap\partial D'$. 
\end{theorem}

In other words, we show that a germ of a biholomorphic mapping between
the boundaries extends analytically to a holomorphic correspondence of
the domains. By the example above, the extension will not be in
general single-valued. Note that we do not require pseudoconvexity of
either $D$ or $D'$. Also $\partial D'$ is not assumed to be simply
connected. If both $D$ and $D'$ have real-algebraic boundary,
i.e. each is globally defined by an irreducible real polynomial, then
we can drop the requirement of simple connectivity of $\partial D$. 

\begin{theorem}\label{t2}
Let $D$ and $D'$ be bounded domains in ${\mathbb C}^n$, $n>1$, with
connected smooth real-algebraic boundaries. Then there exists a proper
holomorphic correspondence between $D$ and $D'$ 
if and only if there exist points $p\in\partial D$,
$p'\in\partial D'$, neighborhoods $U\ni p$ and $U'\ni p'$, and a
biholomorphic map $f:U\to U'$, such that  $f(U\cap\partial D)=
U'\cap\partial D'$. 
\end{theorem}

The proof of Theorem \ref{t2} is simpler than that of Theorem \ref{t1}
due to the fact that when both domains are algebraic, Webster's
theorem \cite{w} can be applied. We give a separate complete proof of
Theorem \ref{t2} in Section 3 to emphasize the ideas of the proof of
Theorem \ref{t1} and the difficulties that arise in the case when one
of the domains is not algebraic.

Local CR-equivalence of the boundaries, which is used in the above
theorems to characterize correspondence equivalence, is a well-studied
subject. Chern and Moser \cite{cm} gave the solution to local
equivalence problem for real-analytic hypersurfaces with
non-degenerate Levi form both in terms of normalization of Taylor
series of the defining functions and in terms of intrinsic
differential-geometric invariants of the hypersurfaces. See also
\cite{se}, \cite{ca} and \cite{t}. Note that for a bounded domain in
${\mathbb C}^n$ with smooth real-analytic boundary, the set of points
where the Levi form is degenerate is a closed nowhere dense set.  Thus
we may reformulate Theorem \ref{t1} in the following way. 

\begin{theorem}
Let $D$ and $D'$ be as in Theorem \ref{t1}. Then $D$ and $D'$ are
correspondence equivalent if and only if there are points
$p\in\partial D$ and $p'\in\partial D'$ such that the Levi form of the
defining functions are non-degenerate at $p$ and $p'$, and the
corresponding Chern-Moser normal forms are equivalent.
\end{theorem}

Theorem \ref{t1} generalizes the result in \cite{p1}, which
states that a bounded, strictly pseudoconvex domain $D\subset{\mathbb
  C}^n$ with connected, simply-connected, real-analytic boundary is
biholomorphically equivalent to  ${\mathbb B}^n={\mathbb B}(0,1)$, the
unit ball centered at the origin, if and only if there
exists a local biholomorphism of the boundaries. It was later
shown in \cite{cj} (see also \cite{ber}) that simple connectivity of 
$\partial D$ can be replaced by simple connectivity of the domain $D$ 
and connectedness of $\partial D$. In \cite{p2} Pinchuk also established
sufficient and necessary conditions of equivalence for
bounded, strictly pseudoconvex domains with 
real-analytic boundaries. Two such domains $D$ and $D'$ are
biholomorphically equivalent if and only if there exist points 
$p\in\partial D$ and $p'\in\partial D'$ such that $\partial D$ and
$\partial D'$ are locally equivalent at $p$ and $p'$. 

If in Theorem \ref{t1} $\partial D$ is not assumed to be simply
connected, then the result is no longer true. Indeed, in a famous
example Burns and Shnider \cite{bs} constructed a domain $\O$ with the
boundary given by 
\begin{equation}\label{bad}
\partial \O=\{z\in{\mathbb C}^2 : \sin(\ln|z_2|) + |z_1|^2=0, \ 
e^{-\pi}\le|z_1|\le 1\},
\end{equation}
which is real-analytic and strictly pseudoconvex but not simply
connected. There exists a mapping $f:{\mathbb B}^2\to \O$ such
that $f$ does not extend even continuously to a point on $\partial
{\mathbb B}^2$. The inverse to $f$ gives a local biholomorphic map
from $\partial \O$ to $\partial {\mathbb B}^2$, but nevertheless
$f^{-1}$ does not extend to a global proper holomorphic correspondence
between $\O$ and ${\mathbb B}^2$. Furthermore, suppose that there exists a
proper correspondence $g:\O\to {\mathbb B}^2$. Since ${\mathbb
B}^2$ is simply connected, $g^{-1}$ is a proper holomorphic mapping,
which extends holomorphically to a neighborhood of $\overline{ {\mathbb
B}^2}$. Let $p\in \partial \O$, $q^1\in f^{-1}(p)$, and $q^2\in
g(p)$. By the result in \cite{bb} $g$ splits at every point in
${\overline \O}$. 
Therefore, with a suitable choice of the branch of $g$ near $p$,  
$\phi:=g\circ f$ defines a local biholomorphic mapping from a
neighborhood $U_1\ni q^1$ to some neighborhood $U_2\ni q^2$. Moreover,
$\phi(U_1\cap \partial {\mathbb B}^2)\subset (U_2\cap \partial
{\mathbb B}^2)$. By the theorem of Poincar\'e-Alexander
(see e.g. \cite{a}), $\phi$ extends to
a global automorphism of ${\mathbb B}^2$. Thus after a biholomorphic
change of coordinates, by the uniqueness theorem $f$ and $g^{-1}$ must
agree in ${\mathbb B}^2$. But $f:{\mathbb B}^2\to \Omega$ is not a
proper map. This contradiction shows that there are no  
proper holomorphic correspondences between $\O$ and ${\mathbb
  B}^2$. Thus the condition of simple connectivity of $\partial D$ in
Theorem \ref{t1} cannot be in general weakened.

One direction in the proof of Theorems \ref{t1} and \ref{t2} is
essentially contained in the work of Berteloot and Sukhov
\cite{bsu}. The proof in the other direction is based on
the idea of extending the mapping $f$ along $\partial D$ as a
holomorphic correspondence.   

It is not known whether Theorem \ref{t1} holds if $\partial D'$ is
real-analytic. The main difficulty is to prove local analytic
continuation of holomorphic mappings along real-analytic hypersurfaces
in the case when hypersurfaces are not strictly pseudoconvex. In
particular, Lemma \ref{l:alongSV} cannot be directly established for
real-analytic hypersurfaces.  

In Section 2 we present background material on Segre varieties and
holomorphic correspondences. In Section 3 we prove a technical
result, important for the proof of the main theorems. Section 4
contains the proof of Theorem~\ref{t2}. In Section 5 we prove local
extendability of holomorphic correspondences along hypersurfaces. Theorem \ref{t1}
is proved in Section 6.

\section{Background Material.}

Let $\G$ be an arbitrary smooth real-analytic hypersurface with a defining 
function $\rho(z,\overline z)$ and let $z^0\in\G$. In a suitable 
neighborhood $U\ni z^0$ to every point $w\in U$ we can associate its 
so-called Segre variety defined as
\begin{equation}
Q_w=\left\{ z\in U: \rho(z,\overline w)=0 \right\},
\end{equation}
where $\rho(z,\overline w)$ is the complexification of the defining
function of $\G$. After a local biholomorphic change of coordinates
near $z^0$, we can find neighborhoods $U_1 \Subset U_2$ of $z^0$, where 
\begin{equation}\label{e:1.00}
U_2={'U_2}\times {{''U_{2}}}\subset {\mathbb C}^{n-1}_{{\ 'z}}\times 
{\mathbb C}_{z_n},
\end{equation}
such that for any $w\in U_1$, $Q_w$ is a closed smooth complex-analytic 
hypersurface in $U_2$.  Here $z=({'z},z_n)$. Furthermore, a Segre variety 
can be written as a graph of a holomorphic function,
\begin{equation}\label{e:implicit}
Q_w=\left\{({'z},z_n)\in {'U_2}\times {''U_{2}} : \ z_n=h({'z},
\overline w)\right\},
\end{equation}
where $h(\cdot,\overline w)$ is holomorphic in ${'U_{2}}$. Following
\cite{dp1} we call $U_1$ and $U_2$ a \textsl{standard} pair of
neighborhoods of $z^0$. A detailed discussion of elementary properties
of Segre varieties can be found in \cite{dw}, \cite{df2}, 
\cite{dp1} or \cite{ber1}.

The map $\lambda:z\to Q_z$ is called the Segre map. We define
$I_w=\l^{-1}\circ \l(w)$. This is equivalent to 
\begin{equation}
I_w=\left\{z\in U_1: Q_z=Q_w\right\}.
\end{equation}
If $\G$ is of finite type in the sense of D'Angelo, or more generally,
if $\G$ is essentially finite, then there exists a neighborhood $U$ of
$\G$ such that for any $w\in U$ the set $I_w$ is finite. Due to the
result in \cite{df1}, this is the case for compact smooth real-analytic
hypersurfaces, in particular for the boundaries of $D$ and $D'$ in
Theorem \ref{t1}. The last observation is crucial for the construction
of proper holomorphic correspondences used throughout this paper. 
We remark that our choice of a standard pair of neighborhoods of
any point $z\in\G$ is always such that for any $w\in U_2$, the set $I_w$ is
finite. 

For the proof of Theorem \ref{t1} we will need the following lemma.  
 
\begin{lemma}\label{l-finite}
Let $\G$ be a compact, smooth, real-algebraic hypersurface. Then there
exist a neighborhood $U$ of $\G$ and an 
integer $m\ge 1$ such that for almost any $z\in U$, $\#I_z = m$. 
\end{lemma} 

\noindent{\it Proof.} Let $P(z,\overline z)$ be the defining
polynomial of $\G$ of degree $d$ in $z$. 
The complexified defining function can be written in the form 
\begin{equation}
P(z,\overline w) = \sum_{|K|=0}^d{a_K(\overline w) z^K},\ \ K=(k_1,\dots, k_n).
\end{equation}
We may consider the projectivized version of the polynomial 
$P:\mathbb P^n \times\mathbb C^n \to \mathbb C$:
\begin{equation}\label{proj-poly}
\tilde P(\tilde\zeta,\overline w)=\zeta_0^d 
\sum_{|K|=0}^d{a_K(\overline w) \left(\frac{\zeta}{\zeta_0}\right)^K}=
\sum_{|K|\le d} a_K(\overline w)\tilde \zeta^K,
\end{equation}
where $z_j=\frac{\zeta_j}{\zeta_0}$, $\zeta=(\zeta_1,\dots,\zeta_n)$,
and $\tilde\zeta = (\zeta_0,\zeta)$. Let $\tilde Q_w = \{\tilde\zeta\in
\mathbb P^n : \tilde P(\tilde\zeta,\overline w)=0\}$. Then $Q_w=\tilde
Q_w\cap \{\zeta_0=1\}$. 

Define the map $\hat \l$ between the set of points in $\mathbb C^n$
and the space of complex hypersurfaces in $\mathbb P^n$, by letting
$\hat \l (w) = \{\tilde\zeta\in\mathbb P^n : \tilde P(\tilde\zeta,
\overline w)=0\}$. Hypersurfaces in $\mathbb P^n$ can be parametrized
by points in $\mathbb P^N$, where $N$ is some large integer, thus
$\hat \l: \mathbb C^n\to \mathbb P^{N}$. Note that each component
of $\hat\l$ is defined by a (antiholomorphic) polynomial. 

It follows that $\hat\l^{-1}\circ\hat\l (w)=I_w$ for every
$w\in\mathbb C^n$, for which $Q_w$ is defined. Indeed, suppose $\xi\in
I_w$. Then 
$\{P(z,\overline\xi)=0\}=\{P(z,\overline w)=0\}$, and therefore, 
\begin{equation}\label{p=p}
\{\tilde P(\tilde\zeta,\overline\xi)=0\}=
\{\tilde P(\tilde\zeta,\overline w)=0\}. 
\end{equation}
Thus $\hat\l(\xi)=\hat\l(w)$. The converse clearly also holds. 

From the reality condition on $P(z,\overline z)$ (see \cite{dw}),
\begin{equation}
I_w\subset \G, {\rm\ \ for\ } w\in\G.
\end{equation}
Since $\G$ is compact, $I_w=\hat\l^{-1}\circ\hat\l (w)$ is finite. Let
$Y=\hat\l(\mathbb C^n)$. Then $\dim Y = n$, and $\hat\l$ 
is a dominating regular map. It follows (see e.g. \cite{m}), that there
exists an algebraic variety $Z\subset Y$ such that for any $q\in Y\setminus Z$,  
\begin{equation}\label{deg}
\#\hat\lambda^{-1} (q)=\deg(\hat\lambda)=m,
\end{equation}
where $m$ is some positive integer. From (\ref{deg}) the assertion
follows. $\square$  

\bigskip

The following statement describes the invariance property of Segre
varieties under holomorphic mappings. It is analogous to the classical 
Schwarz reflection principle in dimension one. 

Suppose that $\G$ and $\G'$ are real-analytic hypersurfaces in
${\mathbb C}^n$, $(U_1, U_2)$ and $(U'_1, U'_2)$  are standard pairs
of neighborhoods for $z_0\in\G$ and $z'_0\in\G'$ respectively. Let $f:
U_2 \to U'_2$ be a holomorphic map, $f(U_1)\subset U'_1$ and $f(\G\cap
U_2)\subset(\G'\cap U'_2)$. Then  
\begin{equation}\label{e-invar}
f(Q_w\cap U_2)\subset Q'_{f(w)}\cap U'_2,\ \ {\rm for\ all\ } w\in U_1.
\end{equation}
Moreover, a similar invariance property also holds for a proper
holomorphic correspondence $f:U_2 \to U'_2$, $f(\G\cap
U_2)\subset(\G'\cap U'_2)$. In this case 
(\ref{e-invar}) indicates that any branch of $f$ maps any point from
$Q_w$ to $Q_{w'}$ for any $w'\in f(w)$. For details, see \cite{dp1} or
\cite{v}. 

Let $f:D\to D'$ be a holomorphic correspondence. We say
that $f$ is {\it irreducible}, if the corresponding analytic set
$A\subset D\times D'$ is irreducible. The condition that $f$ is proper
is equivalent to the condition that 
\begin{equation}
\sup \{ {\rm dist}(f(z),\partial D') \}\to 0,\ \  {\rm as}\ \ 
{\rm dist}(z,\partial D)\to 0.
\end{equation}

Recall that if $A\subset D\times D'$ is a proper holomorphic
correspondence, then $\pi: A\to D$ and $\pi':A\to D'$ are proper. There
exists a complex subvariety $S\subset D$ and a number $m$ such that
\begin{equation} 
f:=\pi'\circ \pi^{-1}=(f^1(z),\dots,f^m(z)), \ \ z\in D,
\end{equation}
where $f^j$ are distinct holomorphic functions in a neighborhood of
$z\in D\setminus S$. The set $S$ is called the {\it branch locus} of
$f$. We say that the correspondence $f$ {\it splits} at
$z\in {\overline D}$ if  there is an open subset $U\ni z$  and
holomorphic maps $f^j:D\cap U\to D'$, $i=1,2,\dots,m$ that
represent $f$.

Given a proper holomorphic correspondence $A$, one can find the system
of canonical defining functions 
\begin{equation}\label{e-canon} 
\Phi_I(z,z') = \sum_{|J|\le m}\phi_{IJ}(z){z'}^J,\ \ |I|=m, \ \ 
(z,z')\in{\mathbb C}^n\times{\mathbb C}^n,
\end{equation}
where $\phi_{IJ}(z)$ are holomorphic on $D$, and $A$ is precisely
the set of common zeros of the functions $\Phi_I(z,z')$. For details
see e.g. \cite{c}.

We define analytic continuation of analytic sets as follows. Let $A$
be a locally complex analytic set in ${\mathbb C}^n$ of pure dimension
$p$. We say that $A$ {\it extends analytically} to an open set 
$U\subset{\mathbb C}^n$ if there
exists a (closed) complex-analytic set $A^*$ in $U$ such that (i)
$\dim A^*=p$, (ii) $A\cap U\subset A^*$ and (iii) every irreducible
component of $A^*$ has a nonempty intersection with $A$ of dimension
$p$. Note that if conditions (i) and (ii) are satisfied, then the last
condition can always be met by removing  certain components of $A^*$.
It follows from the uniqueness theorem that such analytic continuation
of $A$ is uniquely defined. From this we define analytic continuation
of holomorphic correspondences:

\begin{definition}
Let $U$ and $ U'$ be open sets in ${\mathbb C}^n$. Let $f:U\to
U'$ be a holomorphic correspondence, and let $A\subset U\times
U'$ be its graph. We say that $f$ extends as a holomorphic
correspondence to an open set $V$, $U\cap V\not=\varnothing$,
if there exists an open set
$V'\subset {\mathbb C}^n$ such that $A$ extends analytically to a
set $A^*\subset V\times V'$ and $\pi:A^*\to V$ is proper. 
\end{definition}

Note that we can always choose $V'={\mathbb C}^n$ in the definition
above. In general, correspondence $g=\pi'\circ\pi^{-1}:V\to V'$,
where $\pi':A^*\to V'$ is the natural projection, may have more
branches in $U\cap V$ than $f$. The following lemma gives a simple criterion
for the extension to have the same number of branches.

\begin{lemma}\label{l-intersection}
Let $A^*\subset V\times{\mathbb C}^n$ be a holomorphic correspondence
which is an analytic extension of the correspondence $A\subset U\times
{\mathbb C}^n$. Suppose that for any $z\in (V\cap U)$,
\begin{equation}\label{pre}
\#\{\pi^{-1}(z)\}=\#\{\pi^{*-1}(z)\},
\end{equation}
where $\pi:A\to U$ and $\pi^*:A^*\to V$. Then $A\cup A^*$ is a
holomorphic correspondence in $(U\cup V)\times {\mathbb C}^n$.
\end{lemma}

\noindent{\it Proof.} We only need to check that $A\cup A^*$ is closed
in $(U\cup V)\times {\mathbb C}^n$. If not, then there exists a
sequence $\{q^j\}\subset A^*$ such that $q^j\to q^0$ as $j\to\infty$, 
$q^0\in U\times {\mathbb C}^n$, and $q^0\not\in A$. Then $q^j\not\in A$ 
for $j$ sufficiently large. Since by the definition of analytic
continuation of correspondences
$A\cap (U\cap V)\times {\mathbb C}^n\subset A^*$, we have
$$
\#\{\pi^{-1}(\pi^*(q^j))\} < \#\{\pi^{*-1}(\pi^*(q^j))\}.
$$
But this contradicts (\ref{pre}). $\square$

\section{Extension along Segre Varieties.}

Before we prove the main results of the paper we need to establish a
technical result of local extension of holomorphic
correspondence along Segre varieties. This will be used on several
occasions in the proof of the main theorems. 

\begin{lemma}\label{l:alongSV}
Let $\G\subset{\mathbb C}^n$ be a smooth, real-analytic, essentially
finite hypersurface, and let $\G'\subset{\mathbb C}^n$ be a smooth,
real-algebraic, essentially finite hypersurface. Let $0\in\G$, and let
$U_1$, $U_2$ be a sufficiently small standard pair of neighborhoods of
the origin. Let 
$f:U\to{\mathbb C}^n$ be a germ of a holomorphic correspondence such
that  $f(U\cap\G)\subset\G'$, where $U$ is some neighborhood of the
origin. Then there exists a neighborhood $V$ of $Q_0\cap U_1$ and an
analytic set $\L\subset V$, $\dim_{\mathbb C} \L\le n-1$, such that $f$
extends to $V\setminus \L$ as a holomorphic correspondence. 
\end{lemma}

\noindent{\it Proof.} In the case when $\G'$ is strictly pseudoconvex
and $f$ is a germ of a biholomorphic mapping, the result was
established in  \cite{s1}. Here we prove the lemma in a more general
situation.   

Lemma \ref{l:alongSV} only makes sense if $U\subset U_1$. We shrink
$U$ and choose $V$ in such a way that for any $w\in V$, the set
$Q_w\cap U$ is non-empty and connected. Note that if $w\in Q_0$, then
$0\in Q_w$ and $Q_w\cap U\not=\varnothing$.  Let $S\subset U$ be the branch
locus of $f$, and let  
\begin{equation}
\S=\{z\in V: (Q_z\cap U)\subset S\}
\end{equation}
Since $\dim_{\mathbb C} S=n-1$ and $\G$ is essentially finite, $\S$
is a finite set. Define 
\begin{equation}\label{e:A1}
A=\left\{ (w,w')\in (V\setminus\S)\times {\mathbb C}^n :
f\left(Q_w\cap U\right)\subset Q'_{w'} \right\}
\end{equation}
We establish the following facts about the set $A$:
\begin{list}{}{}
\item[(i)] $A$ is not empty
\item[(ii)] $A$ is locally complex analytic
\item[(iii)] $A$ is closed
\item[(iv)] $\S\times{\mathbb C}^n$ is a removable singularity for $A$.
\end{list}

(i) $A\not=\varnothing$ since by the invariance property of Segre
varieties, $A$ contains the graph of $f$.  

\smallskip

(ii) Let $(w,w')\in A$. Consider
an open simply connected set $\O\in(U\setminus S)$ such that
$Q_w\cap\Omega\not=\varnothing$. Then the branches of $f$ are
correctly defined in $\O$. Since $Q_w\cap U$ is connected, the inclusion
$f(Q_w\cap U)\subset Q'_{w'}$ is equivalent to
\begin{equation}\label{e:omega1}
f^j(Q_w\cap\Omega)\subset Q'_{w'}, \ j=1,\dots,m,
\end{equation}
where $f^j$ denote the branches of $f$ in $\Omega$. Note that such
neighborhood $\Omega$ exists for any $w\in V\setminus\S$. 
Inclusion (\ref{e:omega1}) can be written as a system of holomorphic equations as
follows. Let $\rho'(z,\overline z)$ be the defining function of $\G'$. Then
\begin{equation}
\rho'(f^j(z), \overline{w'})=0, \ {\rm\ for\ any}\  z\in (Q_w\cap
\Omega), \ \ j=1,2,\dots,m.
\end{equation}
We can choose $\Omega$ in the form
\begin{equation}
\O={'\O}\times{\O_n}\subset{\mathbb C}^{n-1}_{'z}\times{\mathbb C}_{z_n}
\end{equation}
Combining this with (\ref{e:implicit}) we obtain
\begin{equation}\label{e-system1}
\rho'(f^j({'z},h('z,\overline w)), \overline{w'})=0
\end{equation}
for any $'z\in {'\Omega}$. Then (\ref{e-system1}) is an infinite system of holomorphic
equations in $(w,w')$ thus defining $A$ as a locally complex analytic
variety in $(V\setminus\S)\times{\mathbb C}^n$.

\smallskip

(iii) Let us show now that $A$ is a closed set. Suppose that $(w^j, {w'}^j)\to(w^0,
{w'}^0)$, as $j\to\infty$, where $(w^j, {w'}^j)\in A$ and $(w^0,{w'}^0)\in
(V\setminus\S)\times{\mathbb C}^n$. Then by the definition of $A$,
$f(Q_{w^j}\cap U)\subset Q'_{{w'}^j}$. Since $Q_{w^j}\to
Q_{w^0}$, and $Q'_{{w'}^j}\to Q'_{w^0}$ as $j\to\infty$, by analyticity
$f(Q_{w^0}\cap U)\subset Q'_{{w'}^0}$, which implies that
$(w^0,{w'}^0)\in A$ and thus $A$ is a closed set. Since $A$ is locally
complex-analytic and closed, it is a complex variety in
$(V\setminus\S)\times{\mathbb C}^n$. We now may restrict considerations
only to the irreducible component of $A$ which coincides with the
graph of $f$ at the origin. Then $\dim A = n$.

\smallskip

(iv) Let us show now that $\overline A$ is a complex variety in $V\times{\mathbb
C}^n$. Let $q\in \S$, then 
\begin{equation}
\overline {A}\cap(\{q\}\times{\mathbb C}^n) \subset \{q\}\times
\{z': f(Q_q)\subset Q'_{z'}\}.
\end{equation}
Notice that if $w'\in f(Q_q)\subset Q'_{z'}$, then $z'\in
Q'_{w'}$. Hence the set  $\{z': f(Q_q)\subset Q'_{z'}\}$ has dimension
at most $2n-2$, and $\overline{A}\cap(\S\times{\mathbb C}^n)$ has
Hausdorff $2n$-measure zero. It follows from Bishop's theorem on
removable singularities of analytic sets (see e.g. \cite{c}) that
$\overline{A}$ is an analytic set in $V\times {\mathbb C}^n$.

\smallskip

Thus from (i) - (iv) we conclude that (\ref{e:A1}) defines a
complex-analytic set in $V\times {\mathbb C}^n$ which we denote again
by $A$. Also we observe that since $\G'$ is algebraic, the system of
holomorphic equations in (\ref{e-system1}) is algebraic in $w'$ and
thus we can define the closure of $A$ in $V\times {\mathbb P}^n$. Let
$\pi:A\to V$ and $\pi': A\to {\mathbb P}^n$ be the natural
projections. Since ${\mathbb P}^n$  is compact, $\pi^{-1}(K)$ is
compact for any compact set $K\subset V$, and thus $\pi$ is
proper. 

This, in particular, implies that $\pi(A)=V$. We let 
$\L_1=\pi({\pi'}^{-1}(H_0))$, where $H_0\subset{\mathbb P}^n$
is the hypersurface at infinity. It is easy to see that $\L_1$ is a 
complex analytic set of dimension at most $n-1$. We also consider
the set 
$\L_2:=\pi\{ (w,w')\in A : \dim_{\mathbb C} \pi^{-1}(w) \ge 1\}$. 
It was shown in \cite{s1} Prop.~3.3, that $\L_2$ is a complex-analytic set of
dimension at most $n-2$. Let $\L=\L_1\cup\L_2$. Then $\pi'\circ
\pi^{-1}|_{V\setminus\L}$ is the desired extension of $f$ as a holomorphic
correspondence. $\square$

\section{Proof of Theorem \ref{t2}.}

For completeness let us repeat the argument of \cite{bsu} to prove the 
``only if'' direction of Theorems \ref{t1} and \ref{t2}. Suppose that
$f:D\to D'$ is a proper holomorphic correspondence. Let us show
that $\partial D$ and $\partial D'$ are locally CR-equivalent. 

If $D$ is not pseudoconvex, then for $p\in\widehat D$, there exists
a neighborhood $U\ni p$ such that all the functions in the
representation (\ref{e-canon}) of $f$ extend holomorphically to
$U$. Here $\widehat D$ refers to the envelope of holomorphy of $D$.
Moreover, we can replace $p$ by a nearby point $q\in
U\cap\partial D$ so that $f$ splits at $q$ and at least one of the
holomorphic mappings of the splitting is biholomorphic at $q$.

If $D$ is pseudoconvex, then $D'$ is also pseudoconvex. By \cite{bsu}
$f$ extends continuously to $\partial D$ and we can choose $p\in
\partial D$ such
that $f$ splits in some neighborhood $U\ni p$ to holomorphic mappings
$f^{j}: D\cap U\to D'$, $j=1,\dots,m$. Since $f^{-1}:D'\to D$
also extends continuously to $\partial D'$, the set $\{f^{-1}\circ
f(p)\}$ is finite. Therefore, by \cite{be3}, $f^j$ extend smoothly to
$\partial D\cap U$. It follows that $f^j$ extend holomorphically to a
neighborhood of $p$ by \cite{be} and \cite{df2}. Finally, choose $q\in
U\cap\partial D$ such that $f^j$ is biholomorphic at $q$ for some $j$.

\medskip

To prove Theorem \ref{t2} in the other direction, consider a
neighborhood $U$ of $p\in\partial D$ and a biholomorphic map
$f:U\to{\mathbb C}^n$ such that $f(U\cap\partial D)\subset \partial
D'$. Let us show that $f$ extends to a proper holomorphic
correspondence $F:D\to D'$.

Let $\G=\partial D$ and $\G'=\partial D'$.
Since the set of Levi non-degenerate points is dense in $\G$, by
Webster's theorem \cite{w}, $f$ extends to an algebraic mapping,
i.e. the graph of $f$ is contained in an algebraic variety
$X\subset{\mathbb C}^n \times {\mathbb C}^n$ of dimension $n$. Without
loss of generality assume that $X$ is irreducible, as otherwise
consider only the irreducible component of $X$ containing $\G_f$, the
graph of the mapping $f$. 

Let $E=\{ z\in \mathbb C^n : \dim \pi^{-1}(z)>0\}$, where $\pi:X\to 
\mathbb C^n$ is the coordinate projection to the first component. Then 
$E$ is an algebraic variety in $\mathbb C^n$. Let $f:{\mathbb C}^n
\setminus E \to {\mathbb C}^n$  
now denote the multiple valued map corresponding to $X$. Let $S\subset{\mathbb
  C}^n\setminus E$ be the branch locus of $f$, in other words, for any
$z\in S$ the coordinate projection onto the first component is not
locally biholomorphic near $\pi^{-1}(z)$. To prove Theorem \ref{t2} it is
enough to show that $E\cap \G=\varnothing$. 

\begin{lemma}\label{l2.1}
Let $p\in\G$. If $Q_p\not\subset E$, then $p\not\in E$.
\end{lemma}

\noindent{\it Proof.}  Suppose, on the contrary, that $p\in E$.
Since $Q_p\not\subset E$, there exist a point $b\in Q_p$ and
a small neighborhood $U_b$ of $b$ such that $U_b\cap E = \varnothing$.
Choose neighborhoods $U_b$ and $U_p$ such that for any $z\in U_p$, the
set $Q_z\cap U_b$ is non-empty and connected. Let
\begin{equation}
\S=\{z\subset  U_p: Q_z\cap U_b\subset S\}.
\end{equation}
Similar to (\ref{e:A1}), consider the set
\begin{equation}
A=\left\{ (w,w')\in (U_p\setminus\S)\times {\mathbb C}^n :
f\left(Q_w\cap U_b\right)\subset Q'_{w'} \right\}.
\end{equation}

Then $A\not=\varnothing$. Indeed, since $\dim_{\mathbb C} E\le n-1$,
there exists a sequence of points 
$\{p^j\}\subset (U_p\cap\G)\setminus (E\cup\S)$ 
such that $p^j\to p$ as $j\to\infty$. By the invariance property
of holomorphic correspondences, for every $p^j$ there exists a
neighborhood $U_j\ni p^j$ such that $f(Q_{p^j}\cap U_j)\subset
Q'_{{p'}^j}$, where ${p'}^j\in f(p^j)$. But this implies that
$f(Q_{p^j}\cap U_b)\subset Q'_{{p'}^j}$, and therefore 
$(p^j,{p'}^j)\in A$. Moreover, it follows that 
\begin{equation}\label{e:A=X}
A|_{U_{j}\times{\mathbb C}^n}=X|_{U_{j}\times{\mathbb C}^n}, \ \
j=1,2,\dots,m. 
\end{equation}

Similar to the proof of Lemma \ref{l:alongSV}, one can show that $A$ is
a complex analytic variety in  $(U_p\setminus\S)\times{\mathbb C}^n$,
and that $\S\times{\mathbb C}^n$ is a removable singularity for
$A\subset U_p\times{\mathbb C}^n$. Denote the closure of $A$
in $U_p\times{\mathbb C}^n$ again by $A$. 

Without loss of generality we assume that $A$ is irreducible,
therefore in view of (\ref{e:A=X}) we conclude that 
$A|_{U_{p}\times{\mathbb C}^n}=X|_{U_{p}\times{\mathbb C}^n}$.
Let $\hat f$ be a multiple valued mapping corresponding to $A$. Then
by analyticity, there exists $p'\in\G'\cap \hat f(p)$. Moreover, by
construction, $\hat f(p) = I'_{p'}$. By \cite{dw},
\begin{equation}\label{e-inG}
I'_{z'}\subset\G', \ \ {\rm  for\  any\  } z'\in\G'.
\end{equation}
Now choose $U_p$ so small
that $\overline A \cap (U_p \times \partial U') =\varnothing$, where
$U'$ is a neighborhood of $\G'$.
This is always possible, since otherwise there exists a
sequence of points $\{(z^j, {z'}^j), j=1,2,\dots\}\subset A$, such that $z^j\to p$ and 
${z'}^j\to {z'}^0 \in\partial U'$ as $j\to\infty$. Then $(p,{z'}^0)\in A$ and
${z'}^0\not\in \G'$. But this contradicts (\ref{e-inG}).

This shows that $\hat f:U_p\to U'$ is a holomorphic correspondence
extending $f$, which contradicts the assumption $p\in E$. $\square$ 

\begin{lemma}\label{l2.2}
Let $p\in\G$. Then there exists a change of variables, which is
biholomorphic near $\overline {D'}$, such that in the new coordinate system
$Q_p\not\subset E$. 
\end{lemma}

{\it Proof.} Suppose that  $Q_p\subset E$. Then we find a point
$a\in(\G\setminus E)$ such that $Q_a\cap Q_p\not=\varnothing$. The
existence of such $a$ follows, for example, from \cite{s1}
Prop~4.1. (Note, that $\dim E\cap \G \le 2n-3$). By Lemma
\ref{l:alongSV} the germ of the correspondence $f$
defined at $a$, extends holomorphically to a neighborhood $V$ of
$Q_a$. Let $\L_1$ and $\L_2$ be as in Lemma~\ref{l:alongSV}. 
Since $\dim \L_2<n-1$, we may assume that
$(Q_p\cap V)\not\subset \L_2$. If $(Q_p\cap V)\subset\L_1$, we can
perform a linear-fractional transformation such that $H_0$ is mapped
onto another complex hyperplane $H\subset{\mathbb P}^n$ and such that
$H\cap\G'=\varnothing$. Note that after such transformation $\G'$
remains compact in ${\mathbb C}^n$. Then we may also assume that
$(Q_p\cap V)\not\subset \L_1$. Thus holomorphic extension along $Q_a$
defines $f$ on a non-empty set in $Q_p$. $\square$ 

\bigskip

Theorem \ref{t2} now follows. Indeed, from Lemmas \ref{l2.1} and
\ref{l2.2} we conclude that $E\cap\G=\varnothing$. Since $D$ is
bounded, $D\cap E=\varnothing$, and $X\cap(D\times D')$ defines a
proper holomorphic correspondence from $D$ to $D'$. 

\section{Local Extension.}

To prove Theorem \ref{t1} we first establish local extension of
holomorphic  correspondences.

\begin{definition}
Let $\G$ and $\G'$ be smooth, real-analytic hypersurfaces in ${\mathbb
 C}^n$. Let $f: U\to \mathbb C^n$ be a holomorphic correspondence such that
$f(U\cap\G)\subset \G'$. Then $f$ is called {\it complete} if for any
 $z\in U\cap\G$, $f(Q_z\cap U)\subset Q'_{z'}$ and $f(z)=I_{z'}$.
\end{definition}

By the invariance property of Segre varieties, $f(_{z}Q_z)\subset
Q_{f(z)}$, where $_{z}Q_z$ is the germ of $Q_z$ at $z$, for any $z\in
U\cap \G$. The condition $f(Q_z\cap U)\subset Q'_{z'}$ in the
definition is somewhat stronger: it states that every connected
component of $Q_z\cap U$ is mapped by $f$ into the same Segre
variety. Note that in general Segre varieties are defined only
locally, while the set $U$ can be relatively large. In this case the
inclusion $f(Q_z\cap U)\subset Q'_{z'}$ should be considered only in a
suitable standard pair of neighborhoods of $z$.

The condition $f(z)=I_{z'}$ in the definition above indicates that $f$
has the maximal possible number of branches. It is convenient to
establish analytic continuation of complete correspondences, as such
continuation does not introduce additional branches. 

\begin{lemma}\label{l3.2}\label{l:sublocal}
Let $f:U\to {\mathbb C}^n$ be a complete holomorphic correspondence,
$f(\G\cap U)\subset \G'$, where $\G=\partial D$ and $\G'=\partial D'$,
$D$ and $D'$ are as in Theorem \ref{t1}. Suppose $p\in \partial U\cap\G$ is
such that $Q_p\cap U\ne\varnothing$. Then there exists a neighborhood
$U_p$ of $p$ such that $f$ extends to a holomorphic correspondence
$\hat f: U_p\to {\mathbb C}^n$.
\end{lemma}

\noindent{\it Proof.} The proof of this lemma repeats that of Lemma \ref{l2.1}.
Let $b\in Q_p\cap U$. Consider a small neighborhood $U_b$ of $b$, $U_b\subset U$,
and a neighborhood $U_p$ of $p$ such that for any $z\in U_b$, the set $Q_z\cap
U_p$ is non-empty and connected. As before, let $S\subset U$ be the branch
locus of $f$, and $\S=\{z\subset  U_p: Q_z\cap U_b\subset S\}$.
Define
\begin{equation}\label{e:A}
A=\left\{ (w,w')\in (U_p\setminus\S)\times {\mathbb C}^n :
f\left(Q_w\cap U_b\right)\subset Q'_{w'} \right\}.
\end{equation}
Observe that since $f$ is complete, for any $w\in U \cap\G$, the inclusion
$f(Q_w\cap U_b)\subset Q'_{w'}$ implies that $f(Q_w\cap U)\subset
Q'_{w'}$. In particular, this holds for any $w$ arbitrary close to $p$.
Therefore $A$ is
well-defined if the neighborhood $U_p$ is chosen sufficiently small. 
Analogously to Lemma \ref{l2.1}, one can show that $A$ is a non-empty
closed complex analytic set in $(U_p\setminus\S)\times{\mathbb C}^n$. 
Similar argument also shows that
$\S\times{\mathbb C}^n$ is a removable singularity for $A$, and thus
$\overline A$ defines a closed-complex analytic set in
$U_p\times{\mathbb C}^n$. 

Let us show now that $A$ defines a holomorphic correspondence
$\hat f:U_p\to U'$, where $U'$ is a suitable neighborhood of
$\G'$. 

Consider the closure of $A$ in $U_p\times {\mathbb P}^n$. Recall, that
since ${\mathbb P}^n$ is compact, the projection $\pi:\overline A \to
U_p$ is proper. In particular, $\pi (\overline A)= U_p$. 
Let $U'$ be a neighborhood of $\G'$ as in Lemma
\ref{l-finite}. To simplify the notation, denote the restriction of
$\overline A$ to $U_p \times U'$ again by $A$. Let $\pi: A\to U_p$ and
$\pi': A\to U'$ be the natural projections, and let $\hat f=\pi'\circ
\pi$. 

Let $Z={\pi'}^{-1}(\G')$. Then since $f(\G\cap U)\subset \G'$,
$\pi^{-1}(\G\cap U\cap U_p)\subset Z$. Therefore there exists at least
one irreducible component of $\pi^{-1}(\G\cap U_p)$ which is contained in
$Z$. Thus for any $z\in\G\cap U_p$, there exists $z'\in\G'$ such that
$z'\in \hat f(z)$. By construction, if $z\in U_p$ and $z'\in \hat
f(z)$, then $\hat f(z) = I'_{z'}$. In view of (\ref{e-inG}) we
conclude that $\hat f(\G\cap U_p)\subset \G'$. Now the same argument
as in Lemma \ref{l2.1} shows that $U_p$ can be chosen so small that
$\hat f$ is a holomorphic correspondence. $\square$

\begin{theorem}[Local extension] \label{t:local}
Let $D$ and $D'$ be as in Theorem \ref{t1}, $\G=\partial D$ and 
$\G'=\partial D'$. Let
$f:U\to {\mathbb C}^n$ be a complete holomorphic correspondence,
such that $f(\G\cap U)\subset \G'$, where $\G\cap U$ is connected, and
$\G\cap \partial U$ is a smooth submanifold. Let $p\in\partial U\cap \G$. Then
there exists a neighborhood $U_p$ of the point $p$ such that $f$
extends to a holomorphic correspondence $\hat f:U_p\to {\mathbb
  C}^n$. Moreover, $\hat f|_{U\cap U_p} = f|_{U\cap U_p}$, and the
resulting correspondence $F:U\cup U_p \to {\mathbb C}^n$ is
complete. 
\end{theorem}

\noindent{\it Proof.}
We call a point $p\in\partial U\cap \G$ {\it regular}, if $\partial
U\cap \G$ is a generic submanifold of $\G$, i.e. 
$T^c_p(\partial U\cap\G)=n-2$. We prove the theorem in three
steps. First we prove the result under the assumption $Q_p\cap
U\not=\varnothing$, then for regular points in $\partial U\cap\G$, and
finally for arbitrary $p\in \partial U\cap\G$.

\smallskip

{\it Step 1.} Suppose that $Q_p\cap U\not=\varnothing$. 
Then by Lemma \ref{l:sublocal}
$f$ extends as a holomorphic correspondence $\hat f$ to some
neighborhood $U_p$ of $p$. It follows from the construction that for any
$z\in U\cap U_p$ the number of preimages of $f(z)$ and $\hat f(z)$ is
the same. Thus by Lemma~\ref{l-intersection}, $f$ and $\hat f$ define a
holomorphic correspondence in $U\cup U_p$. Denote this correspondence
by $F$.

We now show that $F$ is also complete in $U\cup U_p$.
Since $f$ is complete, for any $z\in U\cap U_p\cap\G$, arbitrarily close to $p$,
$f(Q_z\cap U)\subset Q'_{f(z)}$. Thus if $U_p$ is chosen sufficiently
small, then for any $z\in U_p\cap\G$,
\begin{equation}
F(Q_z\cap(U\cup U_p))\subset Q'_{F(z)}. 
\end{equation}
Suppose now that there exists some
point $z$ in $(U\cup U_p)\cap\G$ such that not all components of
$Q_z\cap (U\cup U_p)$ are mapped by $F$ into the same Segre
variety. From the argument above, $z\notin U_p$. Since $U\cap\G$ is
connected, there exists a simple smooth curve $\gamma\subset \G\cap U$ 
connecting $z$ and $p$. By Lemma \ref{l:alongSV} for every point
$\zeta\in\g$, the germ of a correspondence $F$ at $\zeta$ extends as a
holomorphic correspondence along the Segre variety $Q_\zeta$. 
Moreover, for $\zeta\in\g$ which are close to $p$, the extension of $F$
along $Q_\zeta$ coincides with the correspondence $f$ in $U$
(even if $Q_\zeta\cap U$ is disconnected). Since
$\cup_{\zeta\in\g}Q_\zeta$ is connected, this property holds for all
$\zeta\in\g$. The extension of $F$ along $Q_\zeta$ clearly maps $Q_\zeta$
into $Q'_{F(\zeta)}$, and therefore $Q_\zeta\cap U$ is mapped by $f$ into the same 
Segre variety. But this contradicts the assumption that the components
of $Q_z\cap U$ are mapped into different Segre varieties. This shows
that $F$ is also a complete correspondence.

\smallskip

{\it Step 2.}
Suppose now that $Q_p\cap U=\varnothing$, but $p$ is a regular point. 
Then by  \cite{s1} Prop. 4.1, there exists a point $a\in U$ such that
$Q_a\cap Q_p\not=\varnothing$. We now apply Lemma \ref{l:alongSV} to
extend the germ of the correspondence at $a$ along $Q_a$. We note that
such extension along $Q_a$ may not in general define a complete
correspondence, since apriori $Q_a\cap \G$ may be disconnected from $U\cap\G$.
Let $\L$ be as in Lemma \ref{l:alongSV}. Then after performing, if
necessary, a linear-fractional transformation in the target space, we
can find a point $b\in Q_p\cap Q_a$, such that $b\notin \L$. Let $U_b$
be a small neighborhood of $b$ such that $U_b\cap\L=\varnothing$ and
$f$ extends to $U_b$ as a holomorphic correspondence $f_b$. Then for any
$z\in U\cap\G$ such that $Q_z\cap U_b\not=\varnothing$, the sets
$f(Q_z\cap U)$ and $f_b(Q_z\cap U_b)$ are contained in the same Segre
variety. Indeed, if not, then we can connect $a$ and $z$ by a smooth
path $\g\subset\G\cap U$ and apply the argument that we used to prove
completeness of $F$ in Step 1. 

Now the same proof as in Step 1 shows that $f$ extends as a
holomorphic correspondence to some neighborhood of $p$, and that the
resulting extension is also complete.

\smallskip

{\it Step 3.} Suppose now that $p\in \partial U\cap \G$ is not a
regular point. Let 
\begin{equation}
M=\left\{z\in\partial U\cap\G: T_z(\partial U\cap\G)=T^c_z(\G)\right\}.
\end{equation}
It is easy to see
that $M$ is a locally real-analytic subset of $\G$. Moreover, since
$\G$ is essentially finite, $\dim M < 2n-2$. Choose the coordinate
system such that $p=0$ and the defining function of $\G$ is given
in the so-called normal form (see \cite{cm}):
\begin{equation}
\rho(z,\overline z)= 2x_n+\sum_{|k|,|l|\ge
  1}\rho_{k,l}(y_n)('z)^k(\overline{'z})^l, 
\end{equation}
where $'z=(z_1,\dots,z_{n-1})$. Since the extendability of $f$
through regular points is already established, after possibly an
additional change of variables,
we may assume that $f$ extends as a holomorphic
correspondence to the points $\{z\in\partial U\cap\G: x_1>0\}$. Let
$L_c$ denote the family of real hyperplanes in the form $\{z\in
{\mathbb C}^n : x_1=c\}$. Then there exists $\epsilon>0$ such for any
$c\in[-\epsilon,\epsilon]$, 
\begin{equation}\label{good-c}
T^c_z(\G) \not= T_z(L_c\cap\G), {\rm \ \ for\ any\ } z\in
L_c\cap\G\cap {\mathbb B}(0,\epsilon).
\end{equation}

Let $\Omega_{c,\delta}$ be the intersection of $\G$, the
$\delta$-neighborhood of $x_1$-axis and the set bounded by $L_c$ and
$L_{c+\delta}$, that is
\begin{equation}
\Omega_{c,\delta}=\left\{
z\in\G\cap {\mathbb B}(0,\epsilon) : 
c<x_1<c+\delta, \ y_1^2+\sum_{j=2}^n|z_j|^2<\delta
\right\}.
\end{equation}

Then there exist $\delta>0$ and $c>0$, such that $f$ extends as a
holomorphic correspondence to a neighborhood  of the set
$\Omega_{c,\delta}$. Since $L_c\cap\G$ consists only of regular
points, from Steps 1 and 2 we conclude that $f$ extends to a
neighborhood of any point in $L_c\cap\G$ that belongs to the boundary
of $\Omega_{c,\delta}$. Let $c_0$ be the smallest number such that $f$ 
extends past $L_c\cap\G$. Then from (\ref{good-c}) and previous steps,
$c_0<0$, and therefore, $f$ extends to a neighborhood of the
origin. $\square$   

\section{Proof of Theorem \ref{t1}.}

The proof of the local equivalence of boundaries $\partial D$ and
$\partial D'$ is equivalent to that of Theorem \ref{t2}.

\medskip

To prove the theorem in the other direction let us first show that a germ
of a biholomorphic map $f:U\to {\mathbb C}^n$, $f(U\cap\G)\subset\G'$
can be replaced by a complete correspondence. Without loss of
generality we assume that $0\in U\cap \G$. We choose a neighborhood
$U_0$ of the origin and shrink $U$ in such a way, that $Q_w\cap U$ is
non-empty and connected for any $w\in U_0$. Define
\begin{equation}\label{comp}
A=\left\{
(w,w')\in U_0\times{\mathbb C}^n: f(Q_w\cap U)\subset Q'_{w'}
\right\}.
\end{equation}
Then (\ref{comp}) defines a holomorphic correspondence, which in particular
contains the germ of the graph of $f$ at the origin. Let $\hat f$ be
the multiple valued mapping corresponding to $A$. Then by construction and
from (\ref{e-inG}), $w'\in\hat f(w)$ implies $\hat f(w)=I'_{w'}$. Thus
$\hat f$ is a complete correspondence.

If $\partial U_0$ is smooth, then By Theorem \ref{t:local} we can
locally extend $\hat f$ along $\G$ past the boundary of
$U_0\cap\G$ to a larger open set $\Omega$. However, local extension in general does
not imply that $A$ is a closed set in $\Omega\times{\mathbb C}^n$.
Indeed, there may exist a point $p\in \G\cap\partial \Omega$ such
that for any sufficiently small neighborhood 
$V$ of $p$, $\Omega\cap V\cap\G$ consists of two connected components, say
$\G_1$ and $\G_2$. Then local extension from $\G_1$ to a neighborhood
of $p$ may not coincide with the correspondence $\hat f$ defined in
$\G_2$. Therefore, local extension past the boundary of $\Omega\cap\G$ 
does not lead to a correspondence defined globally in a neighborhood of
$\overline\Omega\times \mathbb C^n$.
Note that this cannot happen if $\Omega$ is sufficiently small.

We now show that $\hat f$ extends analytically along any path on $\G$. 

\begin{lemma}\label{along-paths}
Let $\g:[0,1]\to\G$ be a simple curve without self-intersections, and
$\g(0)=0$. Then there exist a neighborhood $V$ of $\g$ and a
holomorphic correspondence $F:V\to\mathbb C^n$ which extends $\hat f$.
\end{lemma}

\noindent{\it Proof.} Suppose that $\hat f$ does not extend along
$\g$. Then let $\zeta\in\gamma$ be the first point to which
$\hat f$ does not extend. Let $\epsilon_0>0$ be so small that $\mathbb
B(\zeta,\epsilon)\cap\G$ is connected and simply connected for any
$\epsilon \le \epsilon_0$. Choose a point $z\in B(\zeta,\epsilon_0/2)\cap
\g$ to which $\hat f$ extends. Let $\delta$ be the largest positive
number such that $\hat f$ extends holomorphically to $\mathbb B(z,\delta)\cap\G$. 
By Theorem \ref{t:local} $\hat f$ extends to a
neighborhood of every point in $\partial\mathbb
B(z,\delta)\cap\G$. Moreover, if $\mathbb B(z,\delta)\subset \mathbb
B(\zeta,\epsilon_0)$, then the extension of $\hat f$ is a closed
complex analytic set. Thus $\delta>\epsilon/2$. This shows that $\hat
f$ also extends to $\zeta$, and therefore extends along $\g$. $\square$
\bigskip

Note that analytic continuation of $\hat f$ along $\g$ in Lemma
\ref{along-paths} always yields a complete correspondence. 

The Monodromy theorem cannot be directly applied for multiple valued 
mappings, and we need to show that analytic continuation is
independent of the choice of a curve connecting two points on $\G$.

\begin{lemma}
Suppose that $\g\subset\G$ is a Jordan curve $\g(0)=\g(1)=0$.
Let $F$ be the holomorphic correspondence defined near the origin and
obtained by analytic continuation of $\hat f$ along $\g$. Then $F=\hat f$ in some
neighborhood of the origin.   
\end{lemma}

\noindent{\it Proof.}
Since $\G$ is simply connected and compact, there exists $\epsilon_0 >0$
such that for any $z\in\G$, $\mathbb B(z,\epsilon)\cap\G$ is connected
and simply connected for any $\epsilon\le\epsilon_0$.

Let $\phi$ be the homotopy map, that is 
$\phi(t,\tau):I\times I\to\G$, $\phi(t,0)=\g(t)$, $\phi(t,1)\equiv 0$,
$I=[0,1]$. Let
$\{(t_j,\tau_k)\in I\times I$, $j,k=0,1,2,\dots,m\}$ be the set of
points satisfying:
 
\begin{enumerate}
\item[(i)] $t_0=\tau_0=0$, $t_m=\tau_m=1$,
\item[(ii)]
$\{\phi(t,\tau_k): t_j\le t \le t_{j+1}\} \subset 
\mathbb B(\phi(t_{j},\tau_{k}),\epsilon_0/2)$, \\	
$\{\phi(t_{j},\tau): \tau_k \le \tau \le \tau_{k+1}\} \subset
\mathbb B(\phi(t_{j},\tau_{k}),\epsilon_0/2)$,
for any $j,k<m$.
\end{enumerate}

Suppose that $f$ is a complete holomorphic correspondence 
defined in a ball $B$ of small radius centered at
$\phi(tj,\tau_k)\in\G$. By Theorem \ref{t:local}, $f$ extends
holomorphically past every boundary point of $\partial B\cap\G$. Since
$B(\phi(t_j,\tau_k),\epsilon_0)$ is connected and simply connected, $f$
extends at least to a ball of radius $\epsilon_0/2$. Consider the closed path 
$\g_{j,k}=\{\phi(t,\tau_k): t_j\le t \le t_{j+1}\}\cup \{\phi(t_{j+1},\tau):
\tau_k\le \tau \le \tau_{k+1}\}\cup
\{\phi(t,\tau_{k+1}): t_j\le t \le t_{j+1}\}\cup \{\phi(t_{j},\tau):
\tau_k\le \tau\le \tau_{k+1}\},$
where the second and fourth pieces are traversed in the opposite
direction. Then $\g_{j,k}$ is entirely contained in $\mathbb
(B(\phi(tj,\tau_k),\epsilon_0/2)$. Therefore, analytic continuation of $f$
along $\g_{j,k}$ defines the same correspondence at $\phi(t_j,\tau_k)$. 

Analytic continuation of $\hat f$ along $\g$ can be reduced to
continuation along paths $\g_{j,k}$. Since continuation
along each path $\g_{j,k}$ does not introduce new branches of $\hat
f$, $F=\hat f$. $\square$
\bigskip

Thus simple connectivity of $\G$ implies that the process of local
extension of $\hat f$ leads to a global extension of $\hat f$ to some
neighborhood of $\G$.  Since $\hat f(\G)\subset \G'$, there exist
neighborhoods $U$ of $\G$ and  $U'$ of $\G'$ such that $\hat f: U\to U'$ is a proper
holomorphic correspondence. Let $A$ be the analytic set corresponding to
$\hat f$. By (\ref{e-canon}) there exist functions
$\phi_{IJ}$ holomorphic in $U$ such that $A$ is determined from the
system $\sum_{|J|\le m}\phi_{IJ}(z){z'}^J=0$. By Hartog's theorem all
$\phi_{IJ}$ extend holomorphically to a neighborhood of $\overline D$,
(recall that $\G=\partial D$). This defines a proper holomorphic
correspondence $f:D\to D'$. $\square$

\begin{small}

\end{small}
\end{document}